\documentclass [leqno,12pt]{article}
\usepackage[utf8]{inputenc} 
\usepackage[T1]{fontenc} 
\usepackage{amsmath}
\usepackage{amsfonts}
\usepackage{amsthm}
\usepackage{graphicx}
\usepackage{epstopdf}
\usepackage{color}

\newlength{\oldparindent}
\newcommand{\cL}{{\mathbb {L}}}%espace L^p

\newcommand{\bpf}{\begin{preuve}}
\newcommand{\epf}{ \end{preuve} \medskip}

\newcommand{\benum}{\begin{enumerate}}
\newcommand{\eenum}{\end{enumerate}}

\newcommand{\bitem}{\begin{itemize}}
\newcommand{\eitem}{\end{itemize}}

\newcommand{\brmq}{\begin{rmq}}
\newcommand{\ermq}{\end{rmq}}

\newcommand{\brmqs}{\begin{rmqs}}
\newcommand{\ermqs}{\end{rmqs}}

\newcommand{\bapp}{\begin{application}}
\newcommand{\eapp}{\end{application}}

\newcommand{\bapps}{\begin{applications}}
\newcommand{\eapps}{\end{applications}}

\newcommand{\bdefi}{\begin{definition}}
\newcommand{\edefi}{\end{definition}}

\newcommand{\beq}{\begin{equation}}
\newcommand{\eeq}{\end{equation}}

\def\bpm{\begin{pmatrix}}
\def\epm{\end{pmatrix}}

\newcommand{\bcas}{\begin{cases}}
\newcommand{\ecas}{\end{cases}}

\newcommand{\bex}{\begin{exemp}}
\newcommand{\eex}{\end{exemp}}

\newcommand{\bexs}{\begin{exemps}}
\newcommand{\eexs}{\end{exemps}}

\newcommand{\bprop}{\begin{proposition}}
\newcommand{\eprop}{\end{proposition}}

\newcommand{\bthm}{\begin{theoreme}}
\newcommand{\ethm}{\end{theoreme}}

\newcommand{\bcor}{\begin{corollaire}}
\newcommand{\ecor}{\end{corollaire}}

\newcommand{\blem}{\begin{lemme}}
\newcommand{\elem}{\end{lemme}}

\newcommand{\beqna}{\begin{eqnarray}}
\newcommand{\eeqna}{\end{eqnarray}}

\newcommand{\beqnas}{\begin{eqnarray*}}
\newcommand{\eeqnas}{\end{eqnarray*}}

\newcommand{\cA}{{\mathcal A}}

\definecolor{green}{rgb}{0,.7,.2}
\definecolor{orange}{rgb}{0.9,.5,0}

\newcommand{\LL}{{\rm L}}%generator
\def\tr{\textmd{trace}\,}

\def\Ric{ {\rm{Ric}}} %parentheses have been corrected
\def\det{{ \rm{det}}}  %parentheses have been corrected
\def\Id{{\rm{Id}}} %parentheses have been corrected

\def\cA{{\mathcal A }}

\def\cD{{\mathcal D}}

\def\cG{{\mathcal  G}}
\def\cH{{\mathcal  H}}

\def\cL{{\mathcal L }}
\def\cM{{\mathcal  M}}

\def\cP{{\mathcal P }}

\def\cT{{\mathcal  T}}

\def\cX{{\mathcal X}}

\def\bbC{{\mathbb{C}}}%complexes

\newcommand{\bbR}{{\mathbb {R}}}%reels
\newtheorem{theoreme}{Theorem}[section]
\newtheorem{lemme}[theoreme]{Lemma}
\newtheorem{definition}[theoreme]{Definition}
\newtheorem{proposition}[theoreme]{Proposition}
\newtheorem{corollaire}[theoreme]{Corollary}

\newenvironment{exemp}{\noindent{\bf Example. --- }}{\par}
\newenvironment{exemps}{\noindent{\bf Examples}\benum}{\eenum\par}
\newtheorem{rmq}[theoreme]{Remark}
\newtheorem{rmqs}[theoreme]{Remarks}

\newenvironment{preuve}{\noindent{\it Proof. --- }}
{\hfill\rule{1.3mm}{2mm}\par} 
\newenvironment{application}{\noindent{\bf Application. --- }}{\par}
\newenvironment{applications}{\noindent{\bf Applications. --- 
}\benum}{\eenum\par}

\theoremstyle{definition}

\author{Dominique Bakry, Olfa Zribi\\ {\small Institut de Mathématiques de Toulouse,} \\ {\small Universit\'e Paul Sabatier,} \\ {\small
118 
route de Narbonne,} \\ {\small 31062 Toulouse,} \\ {\small FRANCE}}
\date{}

\makeatletter \renewcommand{\@oddfoot}{\sl \small
 \hfil \thepage\hfil }
\renewcommand{\@oddhead}{\sl \small
 \hfil }
 \makeatother

\title{Curvature dimension bounds on the deltoid model
}

\begin{document}

\maketitle

\def\bZ{\overline{Z}}
\def\bz{\overline{z}}
\def\LL{{\bf\cL}}
\maketitle
\abstract{ The deltoid curve in $\bbR^2$ is the boundary of a domain on which there exist probability measures and orthogonal polynomials for theses measures which are eigenvectors  of diffusion operators.  As such, they may be considered as a two dimensional extension of the classical Jacobi operators. They belong to one of the 11 families of such  bounded domains in $\bbR^2$.  We study the curvature-dimension inequalities associated to these operators, and  deduce various   bounds on the associated polynomials, together with  Sobolev inequalities related to the associated Dirichlet forms}

\section{Introduction}
The deltoid curve (also called Steiner's hypocycloid), see figure page~\pageref{fig:Deltoide}, is one of the bounded domains in $\bbR^2$ on which there exist a probability measure $\mu$ and a symmetric diffusion process, the eigenvectors of it  being  orthogonal polynomials for $\mu$. These orthogonal polynomials have been introduced in~\cite{Koorn1,Koorn2} and appear in the classification of \cite{KrallS}. They appear in \cite{BOZ2013} as one of the eleven models in dimension two for which such polynomials exist.   Moreover, it seems that  is one of the most difficult models to analyse, since there does not exist many geometric interpretation for it. Beyond this, it is also interesting since it is deeply linked with the  analysis  of the $A_2$ root system and of the spectral analysis   of $SU(3)$ matrices.

This deltoid model  and the associated generators provide an interesting object  to check various properties of diffusion operators, since one knows  explicitly the eigenvalues, and has many informations on the eigenvectors. For example,  they have satisfy  recurrence formulas which allows for explicit computations, and in some cases generating functions, see~\cite{Zribi2013}.

Since the associated measures and operators depend  on a real parameter $\lambda>0$ (see equation~\eqref{eq.deltoid}), one may try to understand how functional inequalities and curvature properties depend on this parameter $\lambda$, and hence on geometric properties of the model. 

It turns out that for the specific cases $\lambda=1$ and $\lambda= 4$, one may produce simple geometric interpretations :  in the first case from the Euclidean Laplace operator through the symmetries of the triangular lattice, in the second case from the Casimir operator on $SU(3)$ acting on spectral measures.  The $SU(3)$ model  provides curvature-dimension inequalities for the generic model for $\lambda\geq 1$.  It is not clear however that these inequalities are optimal. It turns out that  they  indeed are. Quite unexpectedly,  the careful investigation of the $CD(\rho,n)$ inequality for this model  does not produce better results than the direct consequence of the  $SU(3)$  inequality. In comparison with the classical case of Jacobi polynomials, which are orthogonal with respect to the measure $C_{a,b}(1-x)^a(1+x)^b dx$ on $(-1,1)$, this situation is similar to the  symmetric case $a=b$, but differs from the dissymmetric case (see~\cite{Bakry96}).

It seems worth to point out at least two interesting features of the computations of curvature-dimension inequalities for this model.  The first aspect concerns the existence of an optimal dimension in the inequality. When one looks for curvature-dimension inequalities on a  compact Riemannian  manifold with dimension $n_0$, but for a reversible  measure which is not the Riemann measure, there is no optimal one. For any $n>n_0$, one may find some constant $\rho(n)$ such that a $CD(\rho(n), n)$ inequality holds. In general, $\rho(n)$ goes to $-\infty$ when $n\to n_0$. It is only for Laplace operators that one may expect some $CD(\rho, n_0)$ inequality. This is not the case here. For any $\lambda>1$, there is a bound $n(\lambda)= 2\lambda$ such that no $CD(\rho, n)$ inequality may hold for $n< n(\lambda)$. However, for this limiting value $n= 2\lambda$, the $CD(\frac{3}{4}(\lambda-1), 2\lambda)$ holds.  Of course, this phenomenon is due to the singularity of the density  of the measure at the boundaries of it's support.

The second aspect concerns the use of appropriate coordinate systems. Since the underlying metric is a flat metric in two dimensions, the curvature-dimension inequality amounts to check for which values $a$ and $b$ one has an inequality of the form 
$$-\nabla\nabla W \geq a \Id + b \nabla W \otimes \nabla W,$$ where $W$ is the logarithm of the density measure with respect to the Riemann measure. 
It turns out that a proper choice of the coordinates leads to very simple formulas, which is not the case if the computation is made through the use of the naive usual coordinates in the Euclidean plane. This comes from the fact that we have indeed at disposal a polynomial structure, expressed through the choice of these coordinates, and the function $W$ satisfies nice relations with respect to this, namely a "boundary equation" described in~\eqref{eq.boundary}. This is a good indication that if one wants to study such inequalities for higher dimensional models, these "polynomial coordinates" should be used instead of the usual ones. 

Then, the curvature-dimension inequalities provide through Sobolev inequalities~\eqref{ineq.sob.gal}  various uniform bounds on the orthogonal polynomials, which turn out, up to some change in the parameters, to be equivalent to the Sobolev inequality itself.
They also provide bounds on kernels constructed from other spectral decompositions, that is for operators which do not necessary commute with our starting operator.

Orthogonal polynomials on the interior of the deltoid curve  belong to the large family of Heckman-Opdam polynomials associated with root systems~\cite{Heckman87, HeckmanOpdam87}, and in the even larger class of MacDonald's polynomials~\cite{Macdo1,Macdo, McDon03}.  As such, they may serve as a guide for more general models of diffusions associated with orthogonal polynomials. One may find some extensive presentation of these pluri-dimensional families of orthogonal polynomials, for example in~\cite{DX,McDon03}. 
 
 The associated  diffusion generator is associated with some reflection group in dimension 2. Operators associated with reflection groups in $\bbR^d$,   known as  Dunkl operators, are extensively studied in the literature, from the points of view of special functions related to Lie group analysis or Hecke algebras, as well as from the point of view of the associated heat equations, or in probability and statistics~\cite{Dunkl,Heck1,Heck2,Rosler,Rosler1}. 
 
 Most of the language and notations related to diffusion operators, and in particular the links between curvature-dimension inequalities and Sobolev inequalities, together with the bounds one may deduce for eigenvectors, are borrowed from~\cite{bglbook}. The fact that Sobolev inequalities are equivalent to bounds on the heat kernel goes back to~\cite{davies89, va-sc-co}, and the relations between curvature dimension inequalities go back to~\cite{bakry-emery85} and are exposed in~\cite{bglbook}, among others. It is classical that  one may deduce from them bounds on the spectral projectors. However, the fact that these bounds may in turn provide Sobolev inequalities, which is the content of Theorem~\ref{recque.sob}, seems quite new, at least to our knowledge, although a similar result concerning logarithmic Sobolev inequalities is exposed in~\cite{bglbook}. Notice that this  recovers a Sobolev inequality with a weaker exponent, not to talk about the optimal constants, which are always out of reach  with this kind of techniques. 
 
 Many of the properties concerning the spectral decomposition of the operator on the deltoid, recurrence formulas for the associated orthogonal polynomials, generating functions, etc., may be found in~\cite{Zribi2013}. We shall not use the results of this paper here, apart the representations coming from symmetry groups in $\bbR^2$ and from $SU(3)$, that we recall briefly in Section~\ref{sec.deltoid}.

The paper is organized as follows. In Section~\ref{sec.cdrn}, we present the general curvature-dimension inequalities and the associated Sobolev inequalities. We show how this provides bounds on the eigenvectors.
Section~\ref{sec.deltoid} is a short introduction to the model associated with the deltoid curve, where we explain the two geometric specific cases. Section~\ref{sec.cd.deltoid} provides the associated curvature-dimension inequalities, first from the $SU(3)$-model, then from direct computation. We provide two approaches for this general case, first in subsection~\ref{subsec.first.proof} using the naive system of coordinates, and then in  subsection~\ref{subsec.simpler.proof} with the adapted system of coordinates.  The choice to present these two approaches, and the striking difference in the complexity of the computation, aims at  underlining the efficiency of the good "polynomial coordinates". Finally, in Section~\ref{sec.bound.on.polyn}, we give  the various bounds on polynomials and operators we are looking for.

\section{Curvature dimension inequalities\label{sec.cdrn}}

We briefly recall in this Section the context of symmetric diffusion operators, following~\cite{BGL}, in a  specific context adapted to our setting.
For a given probability  space $(X, \cX, \mu)$, we suppose given an algebra $\cA$ of functions  such that 
$\cA\subset \cap_{1\leq p<\infty} \cL^p(\mu)$, $\cA$ is dense in $\cL^2(\mu)$, and  which is stable under composition with smooth functions $\Phi$. In our case, $\cA$ may be chosen as the class of restrictions to the domain of smooth functions defined in a neighborhood of it, without any boundary condition. This particular choice for $\cA$ is made possible thanks to  a special property of the operator, which satisfies a "boundary equation", which is our model has the specific form~\eqref{eq.boundary}. It is valid as soon one deals with  operators having  polynomial  eigenvectors on a bounded domain, see~\cite{BOZ}.  In most of the cases, we may as well restrict our attention to polynomials, although this would not be appropriate for Sobolev inequalities.  A bilinear application   
$\Gamma :~\cA\times \cA\mapsto \cA$ is given such that, $\forall f\in \cA$,  $\Gamma(f,f)\geq 0$,  which  satisfies 
$\Gamma(\Phi(f_1, \cdots, f_k), g)= \sum_i \partial_i \Phi \Gamma(f_i,g)$ for any smooth function $\Phi$. A linear operator $\LL$ is defined through
\beq\label{ipp}\int_X f\LL(g) d\mu= -\int_X \Gamma(f,g) d\mu\eeq and we assume that $\LL$ maps $\cA$ into $\cA$.  In this context, all the properties of the model are described by $\Gamma$ and $\mu$, and the model is then entirely described by the triple $(X, \Gamma, \mu)$ (see~\cite{bglbook}). We then extend $\LL$ into a self adjoint operator and we suppose  that $\cA$ is dense in  the domain of $\LL$.

Then, for $f= (f_1, \cdots, f_k)$  and for any smooth function $\Phi$
\beq\label{diff.prop}\LL(\Phi(f))= \sum_1^k \partial_i \Phi(f) \LL(f_i) + \sum_{i,j=1}^k \partial^2_{ij} \Phi(f) \Gamma(f_i,f_j).\eeq
We have 
$$\Gamma(f,g)= \frac{1}{2} \Big(\LL(fg)-f\LL(g)-g\LL(f)\Big).$$

We moreover define the $\Gamma_2$ operator as 
\beq\label{def.gamma2}\Gamma_2(f,g)= \frac{1}{2}\Big( \LL\Gamma(f,g)-\Gamma(f, \LL g)-\Gamma(g, \LL f)\Big).\eeq

Then,  for any parameters $\rho\in \bbR$ and $n\in [1, \infty]$, we say that $\LL$ satisfies a curvature-dimension inequality $CD(\rho,n)$ if and only if
$$\forall f\in \cA, \Gamma_2(f,f) \geq \rho \Gamma(f,f)+ \frac{1}{n} (\LL f)^2.$$

It is worth to observe that the $CD(\rho,n)$ inequality is local.   For a general  elliptic operator on a smooth  manifold $M$ with dimension $n_0$, one may always decompose $\LL$ into $\Delta_g+ \nabla \log V$, where $\Delta_g$ is the Laplace operator associated with the co-metric $(g)$ and $V$ is the density of  $\mu$ with respect to the Riemann measure. In which case, the operator $\Gamma_2$ may be decomposed as 
\beq \label{gamma2} \Gamma_2(f,f) = | \nabla\nabla f|^2 + (\Ric_g-\nabla\nabla \log V)(\nabla f, \nabla f),\eeq where  $\Ric_g$ denotes the Ricci curvature computed for the Riemannian metric associated with $g$,  $\nabla\nabla \log V$ is the Hessian of $\log V$, also computed  in this metric, and $|\nabla\nabla f|^2$ is the Hilbert-Schmidt norm of the Hessian of $f$.

In this case,  the $CD(\rho, n)$ inequality holds if and only  if  $n\geq n_0$ and, when $V$ is not constant, when $n>n_0$ and 
\beq\label{cdrn.on.manifolds}\Ric_g - \nabla\nabla \log V \geq \rho g + \frac{1}{n-n_0} \nabla \log V \otimes \nabla \log V.\eeq

Of course, when $\LL= \Delta_g$, this amounts to $n\geq n_0$ and $\Ric_g \geq \rho g$. In this case, there exists a best choice for both $\rho$ and $n$, namely  for $n$ the dimension of the manifold and for $\rho$ a lower bound on the Ricci tensor, that is the infimum over $\cM$ of the lowest eigenvalue of this tensor. 

In this paper, we are  mainly mainly interested the case where  $\LL=\Delta+\nabla \log(V)$, where $\Delta$ is the Euclidean Laplace operator in some open set of $\bbR^{n_0}$.  In which case the measure $\mu$ is $V dx$, and the  $CD(\rho, n)$ inequality holds if and only if $n \geq n_0$ and 
$$ -\nabla \nabla \log(V) \geq \rho  +\frac{1}{n-n_0} \nabla \log(V) \otimes \nabla \log(V).$$

In order for it to be satisfied, we may look at local inequalities $CD(\rho(x), n(x))$ and try to find such a pair $(\rho(x), n(x))$ for which $\rho(x)$ is bounded below and $n(x)$ bounded above. 
For a generic function $V$, there is no "best"  local inequality in general. The $CD(\rho,n)$ inequality  requires that the symmetric tensor $-\nabla\nabla \log V$ is bounded below by $\rho \Id$. If $\rho_0$ is the best real number such that $-\nabla\nabla \log V\geq \rho_0 \Id$ (that is $\rho_0$ is the lowest eigenvalue at the point $x$ of $-\nabla\nabla \log V$), then the inequality holds as soon as $\rho\geq \rho_0$ and  $(\rho-\rho_0)(n-n_0) \geq |\nabla \log V|^2$. 

But in our case, as already mentioned in the introduction, we are not in this situation. There is a lower bound on the admissible dimension, which is strictly bigger than $n_0$. To understand this phenomenon, one may analyse a bit further this $CD(\rho,n)$ inequality at a given point on the manifold. 

It may happen that, at some point $x$,  the eigenvector   of $-\nabla\nabla \log V$  corresponding to some eigenvalue $\rho_1(x)> \rho_0(x)$   is parallel to $\nabla \log V$.   Let then   in  this case, there is a best choice for both $\rho(x)$ and $n(x)$, which is 
$$\bcas \rho(x)= \rho_0(x),\\ n(x)= n_0+ \frac{1}{\rho_1(x)-\rho_0(x)}|\nabla\log  V|^2.\ecas$$
In the model that we shall consider later, we shall see that this happens  asymptotically on the boundary of the set we are working on, and the constants $n$ and $\rho$ computed at this boundary points are valid for all other points $x$.

When some $CD(\rho,n)$ inequality holds, with  $\rho$ and $n$ constant, and whenever $\rho>0$, and  $2<n<\infty$, then $(X,\Gamma, \mu)$ satisfies a tight Sobolev inequality. For $p= \frac{2n}{n-2}$, and for any $f\in \cA$, we have
\beq\label{ineq.sob.gal}\Big(\int_X f^p d\mu\Big)^{2/p} \leq \int_X f^2 d\mu + \frac{4}{n(n-2)}\frac{n-1}{\rho} \int_X \Gamma(f,f) d\mu.\eeq

More generally, an $n$- dimensional Sobolev inequality is an inequality of the form
\beq\label{def.sob}\|f\|_{2n/(n-2)}^2 \leq A\|f\|_2^2 + C\int \Gamma(f,f)\, d\mu.\eeq
When  $\mu$ is  a probability measure, we say that the inequality is tight when the constant $A$ is $1$, and provided some Sobolev inequality holds, tightness is equivalent to the fact that a Poincaré inequality occurs, which is automatic in our case since the spectrum is discrete (see~\cite{bglbook}). 

When  a Sobolev inequality~\eqref{def.sob} holds, then the associated semigroup $P_t= \exp(t\LL)$ is ultracontractive, that is, for any $q\in [2, \infty]$
\beq\label{ultrac.2q} \|P_t f\|_q\leq \frac{C_1}{t^{\frac{n}{2}(\frac{1}{2}-\frac{1}{q})}} \|f\|_2, ~0<t\leq 1,\eeq
with
$$C_1= \Big(\frac{Cn}{2}+ 2A\Big)^{n/2}.$$  This last constant $C_1$ is not sharp however. 
The bound is valid for $q= \infty$ and indeed, the result for any $q\in (2, \infty)$ is a consequence of the case $q= \infty$ through an interpolation argument.
It turns out that his last ultra contractive bound (for some given $q$, but for any $t\in (0,1)$) is in turn equivalent to the Sobolev inequality.

There is another equivalent version
\beq\label{ultrac.1infty} \|P_t f\|_\infty\leq Ct^{-n/2} \|f\|_1, ~0<t\leq 1,\eeq for which one deduces immediately that the semigroup $P_t$ has a density which is bounded above by $Ct^{-n/2}$.

When   $1\leq n\leq 2$, one may replace Sobolev inequalities by Nash inequalities, which play the same rôle, see remark~\ref{rmq.nash}. However, the best constants that one may deduce from curvature-dimension inequalities for Nash inequalities are not known (see~\cite{BGL}).

 As a consequence of ultracontractive bounds, whenever $f$ is an eigenvector for $\LL$ with eigenvalue $-\lambda$, and provided that $\int f^2 d\mu=1$, one has 
\beq\label{bnd.eigenv.1}\|f\|_q \leq C_1\inf_{t>0} \frac{e^{\lambda t} }{t^{\frac{n}{2}(\frac{1}{2}-\frac{1}{q})}}= C_1 C_{n,q} \lambda^{\frac{n}{2}(\frac{1}{2}-\frac{1}{q})} ,\eeq
with $C_{n,q} = \inf_{s>0} e^{s}s^{-\frac{n}{2}(\frac{1}{2}-\frac{1}{q})}$, 
which follows immediately from the fact that $P_t f= e^{-\lambda t} f$. This applies in particular for $q=\infty$, and produces uniform bounds on the eigenvectors from the knowledge of their $\cL^2(\mu)$ norms.

\section{Diffusion processes on the interior of the deltoid curve\label{sec.deltoid}}

We describe first the operator associated with the deltoid curve associated with a family of orthogonal polynomials. Most of the details may be found in~\cite{Zribi2013}. The deltoid curve is a degree 4 algebraic plane curve which may be parametrized as 
\begin{equation*}
x(t) = 2\cos  t + \cos 2t, \quad y(t) = 2\sin t - \sin 2t
\end{equation*}

\begin{figure}[ht]
 \centering		\includegraphics[width=.5\linewidth]{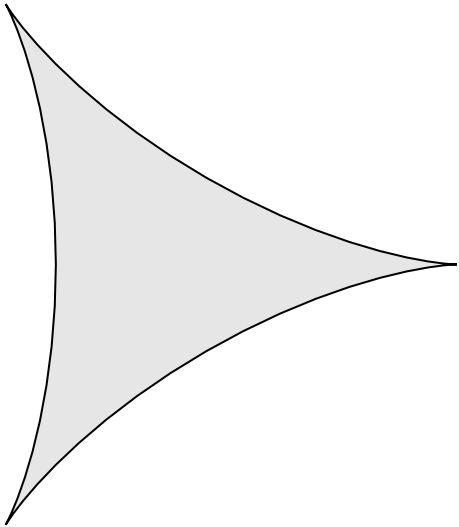}
		\caption{The deltoid domain.}
		\label{fig:Deltoide}
\end{figure}

The connected component  of the complementary of the curve which contains $0$  is a bounded open set, that we refer to as the deltoid domain $\cD$. Indeed, in what follows, we shall work on this domain scaled by the factor $1/3$, which will produce much more convenient formulas.  It turns out that there exist on this domain a one parameter family $\LL^{(\lambda)}$ of symmetric diffusion operator which may be diagonalized in a basis of orthogonal polynomials. It was introduced in~\cite{Koorn1, Koorn2}, and further studied in~\cite{Zribi2013}. This is one of the 11  families of sets in dimension 2  carrying such  diffusion operators, as described in~\cite{BOZ2013}.

In order to describe the operator, and thanks to the diffusion property~\eqref{diff.prop}, it is enough to describe $\Gamma(x,x)$, $\Gamma(x,y)$, $\Gamma(y,y)$, $\LL^{(\lambda)}(x)$ and $\LL^{(\lambda)}(y)$ (the $\Gamma$ operator does not depend on $\lambda$ here).

The symmetric matrix $\bpm \Gamma(x,x)&\Gamma(x,y)\\\Gamma(y,x)&\Gamma(y,y)\epm$ is referred to in what follows as the metric associated with the operator, although properly speaking it is in fact a co-metric.
It is indeed easier to  use  the complex structure of $\bbR^2\simeq \bbC$, and use the complex variables $Z= x+iy$, $\bar Z= x-iy$, with 
$$\LL(Z)= \L(x)+i\LL (y), ~\LL (\bar Z)= \LL (x)-i\LL (y),$$
 and
$$ \Gamma(Z, Z)= \Gamma(x,x)-\Gamma(y,y)+2i\Gamma(x,y), ~\Gamma(Z, \bar Z)=  \Gamma(x,y)+ \Gamma(y,y).$$
         The formulas are much simpler with these variables, 
and $\LL ^{(\lambda)}$ is then described as
\beq\label{eq.deltoid}
\bcas \Gamma(Z,Z) = \bar Z-Z^2,\\ \Gamma(\bar Z, \bar Z)= Z-\bar Z^2, \\ \Gamma(\bar Z, Z)= 1/2(1-Z\bar Z),\\
\LL(Z)= -\lambda Z, \LL(\bar Z)= -\lambda \bar Z,
\ecas
\eeq 
where $\lambda>0$ is a real parameter.

The boundary of this domain turns out to the curve with equation 
\beq\label{defP}P(Z,\bar Z)= \Gamma(Z, \bar Z)^2-\Gamma(Z,Z)\Gamma(\bar Z, \bar Z)=0,\eeq and inside the domain, the associated metric is positive definite, so that it corresponds to some elliptic operator on it. The reversible measure associated with it has density $C_\alpha P(Z, \bar Z)^\alpha$ with respect to the Lebsegue measure, where $\lambda = (6\alpha+5)/2$, and is a probability measure exactly when $\lambda>0$ (we refer to~\cite{Zribi2013} for more details).

There are two particular cases which are worth understanding, namely $\lambda= 1$ and $\lambda= 4$, corresponding to the parameters $\alpha= \pm1/2$. We briefly present those two models, referring to~\cite{Zribi2013} for more details, since we shall make a strong use of them in what follows.

In the first case, one sees that this operator is nothing else that the image of the Euclidean Laplace operator on $\bbR^2$ acting on the functions which are invariant under the symmetries around the lines of a regular triangular lattice. 

For the first one, one considers the three unit  third roots of identity in $\bbC$, say $(e_1,e_2,e_3)= (1,j, \bar j)$. Then, consider the functions  $Z$ and  $z_k: \bbC\mapsto \bbC$ which are defined  as 
\beq\label{def.zk}z_k(z) = e^{i\Re(z\bar e_k)}, ~Z= \frac{1}{3}(z_1+z_2+z_3)
\eeq
They satisfy $|z_k|=1$ and $z_1z_2z_3=1$.

 It is easily seen that, for the Euclidean Laplace operator on $\bbR^2$, $Z$ and $\bar Z$ satisfy the  relations~\eqref{eq.deltoid} with $\lambda=1$. Moreover, the function $Z: \bbC\mapsto \bbC$ is a diffeomorphism between the interior $\cT$ of the triangle $T$ and the deltoid domain $\cD$, where $T$ is one of  the equilateral triangles with containing  the two edges $0$ and $4\pi/3$. The functions which are invariant under the symmetries of the triangular lattice generated by this triangle $T$ are exactly functions of $Z$.
Therefore, the image of $\LL ^{(1)}$ through  $Z^{-1} : \cD\mapsto \cT$ is nothing else that the Laplace operator on $\cT$. We may as well look at the image of the operator $\LL^{(\lambda)}$ and it is then simply
\beq \label{repLtriangle}\LL ^{(\lambda)} = \Delta(f)+ (\alpha+1/2)\nabla\log W\nabla f=\Delta(f)+ \frac{\lambda-1}{3}\nabla\log W\nabla f \eeq 
where $\Delta$ is the usual Laplace operator in $\bbR^2$ and  the function $W$ is defined form the functions $z_j$ described in equation~\eqref{def.zk} as
\beq \label{rep.measure}W= -(z_1-z_2)^2(z_2-z_3)^2(z_3-z_1)^2.\eeq

One should be aware here that thanks to the properties of the functions $z_j$, $-(z_1-z_2)^2(z_2-z_3)^2(z_3-z_1)^2$ is a real valued function taking values in $(0,\infty)$ (and vanishes only at the boundaries of $\cT$). 

This representation provides a way of computing $CD(\rho,n)$ inequalities for $\LL ^{(\lambda)}$, following the description of Section~\ref{sec.cdrn}. 

The second description follows from the Casimir operator on $SU(3)$. This latter group is a semi-simple compact Lie group, and as such as a canonical Laplace (Casimir) operator which commutes (both  from left and  right) to the group action~\cite{Helgason2001,faraut2005} . Namely, in any such semi simple compact Lie group $G$, one considers it's Lie algebra $\cG$, naturally endowed with a Lie algebra structure $\cG\times \cG\mapsto \cG$, $(X,Y)\mapsto [X,Y]$. The Lie algebra structure provides on $\cG$ a natural quadratic form $K$ (the Killing form) as follows : for any element $X\in \cG$, one considers $ad(X) : \cG\mapsto \cG$,  $Y\mapsto [X,Y]$, and 
$K(X,Y) = -\tr\Big(ad(X)ad(Y)\Big)$. It turns out that this quadratic form is positive definite exactly when the group is compact and semi-simple. 
If one considers, for this Killing form,  any orthonormal basis $(X_i)$ in the Lie algebra, the quantity $\sum_i X_i^2$, computed in the enveloping algebra, does not depend on the choice of the basis, and  commutes with any element on the Lie algebra itself (this means that this commutation property depends only on the Lie-algebra structure and not on the way that the elements of this Lie algebra are effectively represented as linear operators). 

Now, to any $X\in\cG$ is associated a first order operator $D_X$ on $G$ defined as follows
\beq\label{vect.field.lie}D_X(f)(g)= \partial_t f(ge^{tX})\mid_{t=0}.\eeq The application $X\mapsto D_X$ is a representation of the Lie algebra into the linear space of vector fields ($[D_X, D_Y]= D_{[X,Y]}$), and any identity in the Lie algebra (on on it's enveloping algebra) translates to an identity on those differential operators.  We work here with the right action $g\mapsto ge^{tX}$ but we could as well work with the left action $g\mapsto e^{tX}g$. For any orthonormal basis $X_i$ of $\cG$ for the Killing form $K$, one defines the Casimir operator 
\beq\label{casimir}\LL = \sum_i D_{X_i}^2.\eeq

It does not depend of the choice of the basis  and commutes with the group action, that is $[\LL, D_X]=0$ for any $X\in \cG$. 

This Killing form  provides an Euclidean quadratic form in the tangent plane at  identity  in $G$  (the Lie algebra $\cG$), which may be translated to the tangent plane  at any point $g\in G$ through the group action, and endows $G$ with a natural Riemanian structure. It turns out that the Casimir operator $\LL $ is also the Laplace operator for this structure. For the group $SU(d)$ that we are interested in, one wants to precisely describe the action of this Casimir operator  on the entries of the matrix $g= (z_{ij})$ in $SU(d)$.  That is, writing the entries $z_{pq}= x_{pq}+ i y_{pq}$, we consider $x_{pq}$ and $y_{pq}$ as functions $G\mapsto \bbR$, and, for any $i,j,k,l$, we want to compute 
$$\LL ^{SU(d)}(x_{ij}), ~\LL ^{SU(d)}(y_{ij}),~ \Gamma^{SU(d)}( x_{ij}, x_{k,l}),~ \Gamma^{SU(d)}( x_{ij}, y_{k,l}),~ \Gamma^{SU(d)}( y_{ij}, x_{k,l}),$$ where $\Gamma^{SU(d)}$ is the square field operator associated with $\LL ^{SU(d)}$. In order to get simpler formulae, it is once again better to work with the complex valued functions $z_{pq}$, writing for such a function $z= x+iy$,
\beqnas&& \LL ^{SU(d)}(z)=\LL ^{SU(d)}(x)+i L^{SU(d)}(y),\\
 &&\Gamma^{SU(d)}(z,z)=\Gamma(x,x)-\Gamma(y,y)+ 2i\Gamma(x,y),\\
  &&\Gamma^{SU(d)}(z, \bar z)= \Gamma(x,x)+ \Gamma(y,y).
 \eeqnas 

 If one denotes by $E_{p,q}$ the matrix with entries  $(E_{p,q})_{i,j}= \delta_{ip}\delta_{jq}$, a 
  base of the Lie algebra of $SU(d)$ is given by 
 $$R_{k,l}=(E_{k,l}-E_{l,k})_{k<l}$$ $$S_{k,l}=i(E_{k,l}+E_{l,k})_{k<l}$$ $$ D_{1,l}=i(E_{1,1}-E_{l,l})_{1<l} .$$
 In order to describe the Casimir operator in a compact form, it is better to work with $D_{k,l}= i(E_{k,k}-E_{l,l})_{k<l}$, and one may write in this way (up to a factor $2$ that will play no rôle in the future)
 \beq\label{def.cas.sud}\LL ^{SU(d)}= \sum_{k<l} (D_{R_{k,l}}^2+D_{S_{k,l}}^2+\frac{2}{n}D_{D_{k,l}}^2).\eeq

 One may compute then the associated vector fields following formula~\eqref{vect.field.lie}, and we get
   $$D_{R_{pq}}=\sum_k z_{kq}\partial_{kp}-z_{kp}\partial_{kq}+  \bar z_{kq}\bar \partial_{kp}-\bar z_{kp}\bar \partial_{kq}$$
$$D_{S_{pq}}= i\Big[\sum_kz_{kq}\partial_{kp}+z_{kp}\partial_{kq}-  \bar z_{kq}\bar \partial_{kp}- \bar z_{kp}\bar \partial_{kq}\Big].$$
$$D_{D_{pq}}= i\Big[\sum_k z_{kp}\partial_{kp}-z_{kq}\partial_{kq}-  \bar z_{kp}\bar \partial_{kp}+\bar z_{kq}\bar \partial_{kq}\Big].$$

With these relations, one may directly compute the action of $\LL ^{SU(d)}$ on the entries of the matrix, and we get
 \beq\label{lapl.sud.coord} \bcas 
 \LL ^{SU(d)}(z_{pq})=\frac{-2(d^2-1)}{d}z_{pq}\\ 
 \LL ^{SU(d)}(\bar z_{pq})=\frac{-2(d^2-1)}{d}\bar z_{pq}\\ 
\Gamma^{SU(d)}(z_{kl},z_{rq})= -2z_{kq}z_{rl}+\frac{2}{d}z_{kl}z_{rq}, \\
\Gamma^{SU(d)}(z_{kl},\bar z_{rq})= 2(\delta_{kr}\delta_{lq} -\frac{1}{d}z_{kl}\bar z_{rq}).
\ecas
\eeq

Now, let us consider the special case $d= 3$ and consider the function $Z: G \mapsto \bbC$ defined by $Z(g)= \frac{1}{3}\tr(g)= \frac{1}{3}(z_{11}+z_{22}+ z_{33})$.  Any function on $G$ depending only on the spectrum of $g$ is a function of $(Z, \bar Z)$, since, denoting by $\lambda_1, \lambda_2, \lambda_3$ the eigenvalues of $g$, with  $|\lambda_i|=1$, $\lambda_1\lambda_2\lambda_3=1$,  such a spectral function is a function of 
$(3Z=\lambda_1+\lambda_2+ \lambda_3, \lambda_1\lambda_2 + \lambda_2\lambda_3+ \lambda_3\lambda_1)$, but  thanks to the properties of the eigenvalues,
$$\lambda_1\lambda_2 + \lambda_2\lambda_3+ \lambda_3\lambda_1= 3\bar Z.$$

Then, the  characteristic polynomial of $g$, that is $P(X)=\det(X-g)$, may be written as 
$$ P(X)=X^3-3Z X^2+ 3\bar Z X -1.$$

It turns out that formulae~\eqref{lapl.sud.coord} produce an easy way to compute the action of $\LL ^{SU(d)}$ on the entries of the characteristic polynomial. For this, we use the following formulae, valid for any square matrix $M$ with entries $(m_{ij})$
\beq\label{der.det.mat} \partial_{m_{ij}}\log \det(M)= M^{ -1}_{ji}, ~\partial^2_{m_{ij}m_{kl} }\log \det(M)= -M^{-1}_{jk}M^{-1}_{li}.\eeq
Using the change of variable formula~\eqref{diff.prop}, when applying equation~\eqref{der.det.mat} to $X\Id -M$  together with  formulae~\eqref{lapl.sud.coord}, one gets, on $SU(d)$ and with $P(X)= \det(X\Id-M)$,
\beqnas\Gamma^{SU(d)}(\log P(X), \log P(Y))=&& (d-\frac{XP'(X)}{P(X)})(d-\frac{YP'(Y)}{P(Y)})\\&& -d\Big(d+ \frac{X+Y}{X-Y}(\frac{YP'(Y)}{P(Y)}\\&&-\frac{XP'(X)}{P(X)})-\frac{ XY}{X-Y}(\frac{P'(Y)}{P(Y)}-\frac{P'(X)}{P(X)})\Big)\\=&&
 XY\Big(\frac{P'(X)}{P(X)}\frac{P'(Y)}{P(Y)}+ \frac{d}{X-Y}(\frac{P'(X)}{P(X)}-\frac{P'(Y)}{P(Y)})\Big).
 \eeqnas
 which in turn gives
 \beq\label{eq.gamm.sud.spect}\Gamma(P(X), P(Y))= XY\Big(P'(X)P'(Y)+ d\frac{P'(X)P(Y)-P'(Y)P(X)}{X-Y}\Big),\eeq
 and
 \beq\label{eq.l.sud.spect}\LL ^{SU(d)}(P)= (1-d^2) X P' + (1+d)X^2 P''.\eeq

One may now compare this for $d=3$ with the formula given by \eqref{eq.deltoid} for $\lambda= 4$, to observe that indeed, with $\LL =\frac{3}{4} \LL ^{SU(3)}$, they give the same result. 

In the end,  we see that $\LL ^{(4)}$ is nothing else than $\frac{3}{4}\LL ^{SU(3)}$ when acting on functions of  $Z= T/3$, where   $T$ is the trace of the matrix.

\section{Curvature dimension for the deltoid model\label{sec.cd.deltoid}}

\subsection{Curvature-dimension inequalities for $SU(3)$.}

It is well known that the Casimir operator of any compact semi simple  Lie group has a constant Ricci curvature (see for example~\cite{gallothulinlafontaine04}, prop. 3.17, or~~\cite{saloff2004convergence}). However, the explicit constant is not straightforward to compute, and for the sake of completeness, we provide for this  an easy way through the use of the $\Gamma_2$  operator. It  relies   on the following observation,  which may be used in other similar situations. Let $E_i$ be some elements of the Lie algebra  and $X_i$ be the associated vector fields. We do not require that the $E_i$ form  an orthonormal basis, since we shall use the representation~\eqref{def.cas.sud} for $SU(3)$ which is not given in such a basis (the elements $D_{ij}$ are not orthogonal and are not linearly independent).

Suppose that $\LL $ is given as  
$$\LL = \sum_i X_i^2, ~ [\LL , X_i]=0, ~\Gamma(f,f)= \sum_i X_i(f)X_i(f).$$
Then, using formula~\eqref{def.gamma2},  one immediately gets 
$$\Gamma_2(f,f)= \sum_{i,j} (X_iX_j f)^2.$$
We decompose $X_iX_j(f)$ into it's symmetric and antisymmetric part 
$$X_iX_j =\frac{1}{2} (X_iX_j f + X_j X_i f)+ \frac{1}{2} [X_i,X_j] (f)= H_{ij} + \frac{1}{2} [X_i,X_j] (f).$$
Then
 $$\Gamma_2(f,f)=\sum_{ij} (H_{ij} f)^2+ \frac{1}{4} \sum_{ij} ([X_i, X_j](f))^2= \sum_{ij} (H_{ij} f)^2+ \frac{1}{2} \sum_{i<j} ([X_i, X_j](f))^2.$$
 
 It turns out that this decomposition coincides exactly with the one provided by~\eqref{gamma2} into the second order part and the first order part of the $\Gamma_2$ tensor. 
 Therefore, on has 
 $$\Ric(f,f) = \frac{1}{2}  \sum_{i<j} ([X_i, X_j](f))^2,$$ and we are bound to compute the commutators of the elements in the basis.
 
 For example, on  $SU(3)$, the commutator table is the following, with $\hat D_{i,j}= \sqrt{2/3}D_{i,j}$, and  $a= \sqrt{2/3}$: 
 
 \begin{tabular}{||c||c|c|c|c|c|c|c|c|c||}
\hline\hline
 0 &$R_{1,2}$&$R_{1,3}$&$R_{2,3}$&$S_{1,2}$&$S_{1,3}$&$S_{2,3}$& $\hat D_{1,2}$&$\hat D_{1,3}$&$\hat D_{2,3}$\\
 \hline\hline
 $R_{1,2}$ & $0$& $-R_{2,3}$&$R_{1,3}$&$2/a\hat D_{1,2}$&$-S_{2,3}$&$S_{1,3}$& $-2aS_{1,2}$&$-aS_{1,2}$&$aS_{1,2}$\\
 \hline
$R_{1,3}$& & 0& $-R_{1,2}$ &$-S_{2,3}$&$2/a\hat D_{1,3}$&$S_{1,2}$&$-aS_{2,3}$&$-2aS_{1,3}$&$-aS_{1,3}$        \\
\hline
$R_{2,3}$& & & 0& $-S_{1,3}$&$S_{1,2}$&$2/a\hat D_{2,3}$&$aS_{2,3}$&  $-aS_{2,3}$ &$-2aS_{2,3}$\\
\hline
$S_{1,2}$ && &&0& $-R_{2,3}$&$-R_{1,3}$ &   $2aR_{1,2}$ &$ aR_{1,2}$ &$-aR_{1,2}$   \\
\hline
$S_{1,3}$    &&& &&0&    $-R_{1,2}$& $aR_{1,3}$& $2aR_{1,3}$ & $aR_{1,3}$     \\
\hline
$S_{2,3}$  &&&& &&0&    $-aR_{2,3}$&    $aR_{2,3}$ &    $2aR_{2,3}$ \\
\hline
$\hat D_{1,2}$ &&&&& &&0&0&0\\
\hline
$\hat D_{1,3}$&&&&&& &&0&0\\
\hline
$\hat D_{2,3}$&&&&&&& &&0\\
\hline\hline
\end{tabular}

We just have to add half the sum of all the squares of the element appearing in the table. We get

   $$(1+3a^2)(R_{2,3}^2+R_{1,3}^2+R_{1,2}^2+S_{1,2}^2+S_{1,3}^2+S_{2,3}^2)+ \frac{2}{a^2}(\hat D_{1,2}^2+\hat D_{2,3}^2+ \hat D_{2,3}^2).$$
 Since $1+3a^2= 2/a^2= 3$, we see that
 the Ricci curvature is constant and it is equal 3, more precisely $\Ric(f,f) = 3 \Gamma(f,f)$.
 
Remembering that the dimension of $SU(3)$ is $8$, we get 

\bprop The Casimir operator $\LL ^{SU(3)}$ defined by equation~\eqref{def.cas.sud} satisfies the optimal $CD(3, 8)$ inequality.

\eprop

From the definition of the $\Gamma_2$ operator and of the curvature-dimension inequality, it is immediate that if if $\LL $ satisfies $CD(\rho,n)$, then $c\LL $ satisfies $CD(c\rho, n)$.  We  therefore see that $\frac{3}{4}\LL ^{SU(3)}$ satisfies a $CD(\frac{9}{4}, 8)$ inequality.

Then, $\LL ^{(4)}$ being  the image of $\frac{3}{4}\LL ^{SU(3)}$ through the map $g\mapsto  Z=\frac{1}{3} (z_{11}+z_{22}+ z_{33})$,  applying the $CD(\rho, n)$ inequality on $SU(3)$ on function of $(Z, \bar Z)$ provided directly the

\bcor  $\LL ^{(4)}$ satisfies the $CD(\frac{9}{4}, 8)$ inequality.

\ecor

\subsection{Curvature dimension for the general deltoid model\label{subsec.first.proof}}

We may now come back to the general deltoid model. Let us write, in the triangle representation,  the decomposition~\eqref{repLtriangle}
$$\LL ^{(\lambda)}= \Delta +\frac{\lambda-1}{3} \nabla \log W,$$ where $W^{(1-\lambda)/3}$ is the density the invariant measure of $\LL ^{(\lambda)}$ with respect to the Lebesgue measure, and $W$ is given in formula~\eqref{rep.measure}s.  

We also know from the general formulation of the $CD(\rho, n)$ inequality~\eqref{cdrn.on.manifolds} that the $CD(\frac{9}{4}, 8)$  for $\lambda= 4$ translates into the following
\beq\label{ineq.gamma2.w}-\nabla\nabla\log W \geq \frac{9}{4} \Id + \frac{1}{6} \nabla\log W \otimes \nabla \log W.\eeq

For $\lambda>1$, multiplying the previous inequality by $\frac{\lambda-1}{3}$,  this in turns gives for the general case provides
\bcor  For any $\lambda\geq 1$, $\LL ^{(\lambda)}$ satisfies a $CD(\frac{3(\lambda-1)}{4}, 2\lambda)$ inequality.
\ecor

Observe  that indeed the limiting case $\lambda=1$ corresponds to the Laplace operator on $\bbR^2$ which satisfies a $CD(0,2)$ inequality.

It is not clear however that this inequality is sharp. Indeed, going from the $CD$ inequality on $SU(3)$ to the same $CD$ inequality for $\LL ^{(4)}$, we may as well have lost information. In general, as we already mentioned, on a smooth compact manifold (with no boundaries) with dimension $n_0$, there is no optimal $CD(\rho, n)$ inequality, and for any $n>n_0$ one may find some $\rho(n_0)$ such that the operator $\Delta + \nabla\log V$ satisfies a  $CD(\rho, n)$ inequality.  Moreover, this does not tell us anything about the case where $\lambda<0$.

As mentioned in the introduction, the optimal computation for this inequality for a generic $\lambda$ is not elementary. In this section, we shall perform directly the computations of $\nabla\nabla \log W$ and $\nabla \log W \otimes \nabla \log W$ to observe that the $SU(3)$ case indeed gives the optimal answer, which is quite surprising. Of course, on $SU(3)$ the $CD(\rho,n)$ inequality is optimal at every point $x\in SU(3)$, while on the projected model, it is optimal only at some point on the boundary of the deltoid domain $\cD$.

We shall  show the following

\bprop ~\label{CD.r.n.brute.force}
\benum
\item
For $\lambda<1$, the operator $\LL ^{(\lambda)}$ does not satisfy any $CD(\rho, \infty)$ inequality.
\item For $\lambda>1$, the operator satisfy no $CD(\rho, n)$ inequality for any $n< 2\lambda$. Moreover, the best constant $\rho$ in the $CD(\rho, 2\lambda)$ inequality is 
$\rho = \frac{3(\lambda-1)}{4}$.
\eenum
\eprop

\bpf 
Everything boils down, for $\lambda >1$, to check inequalities of the form 
$$-\nabla\nabla \log W \geq c_1 \Id + c_2 \nabla \log W \otimes \nabla \log W,$$ and for $\lambda < 1$ to check the the tensor 
$\nabla\nabla \log W$ is not bounded below.

To perform the computations, we shall use the triangle model,  that is move back everything of $\cT$ through the map $Z^{-1}$, since on $\cT$   the metric is the identity and the Hessian is computed in the usual way.

In what follows, we shall use the functions $\log (z)$ for a complex variable $z \neq 0$ without any  precaution about which determination of the argument we chose for the logarithm, since indeed we are only concerned with the one form $d\log (z)$ and it's derivative.

Let us recall that $W= -(z_1-z_2)^2(z_2-z_3)^2(z_3-z_1)^2$, with $z_k= e^{i(E_k.z)}$, where 
  \beq\label{defEi}E_{1}=\begin{pmatrix}
 1\\
0
 \end {pmatrix},
 E_{2}=\begin{pmatrix}
 -\frac{1}{2}\\
\frac{\sqrt3}{2}
 \end {pmatrix},
E_{3}=\begin{pmatrix}
 -\frac{1}{2}\\
-\frac{\sqrt3}{2}
 \end {pmatrix}.\eeq

Then, setting $\sigma= \log(z_1-z_2)+\log(z_2-z_3)+\log(z_3-z_1)$,  $\log W =  2\sigma$, up to some additive (eventually complex) constant,  and we are looking for
\beq\label{ineg.gamma2.local}- \nabla \nabla \sigma  \geq b I_{2}+ a\nabla \sigma   \otimes \nabla \sigma  . \eeq    

What we want to show is first that $\nabla\nabla \sigma$ is not bounded above, and then that the former inequality~\eqref{ineg.gamma2.local} may not hold if $a>1/3$. Moreover, we want to check that for $a= 1/3$, the best lower bound for $b$ is $b= 9/8$. 
It turns out that this is quite technical.

We shall need a few intermediate steps to check this inequality.

\blem\label{lem.derz} We have

 $$\nabla z_{k}=iz_{k}E_{k}, ~
 \nabla \log(z_{p}-z_{q})=i\frac{z_{p}E_{p}-z_{q}E_{q}}{z_{p}-z_{q}},$$
 $$ \nabla\nabla \log(z_p-z_q)=\frac{z_{p}z_{q}}{(z_{p}-z_{q})^2}(E_{p}-E_{q})^{\odot^2}.$$
 \elem
\bpf (Of Lemma~\ref{lem.derz})

The two first identities are immediate. For the last one, we write
\beqnas
\nabla \nabla\log(z_p-z_q)&=&i\nabla(\frac{1}{z_p-z_q}(z_pE_p-z_qE_q))\\
&=& -\frac{1}{z_p-z_q}(z_pE_p^{{\otimes}^2}-z_qE_q^{{\otimes}^2})+\frac{1}{(z_p-z_q)^2}(z_pE_p-z_qE_q)^{{\otimes}^2}\\
\eeqnas

\epf
As a consequence, one has
\bcor\label{cor.eqwder}
\beq\label{eqwder1}\nabla\nabla \log \sigma =[\frac{z_1z_2}{(z_1-z_2)^2}V_1^{\odot^2}+\frac{z_2z_3}{(z_2-z_3)^2}V_2^{\odot^2}+\frac{z_1z_3}{(z_3-z_1)^2}V_3^{\odot^2}].\eeq
\beq\label{eqwder2}\nabla \log\sigma  \otimes \nabla \log \sigma =-[\frac{z_1}{z_1-z_2}V_1+\frac{z_2}{z_2-z_3}V_2+\frac{z_3}{z_3-z_1}V_3]^{\odot^2}.\eeq

where  $$V_1=E_1-E_2, V_2=E_2-E_3, V_3=E_3-E_1$$ and  $E_i,i=1,2,3$ are defined in \eqref{defEi}. \ecor

\brmq Observe that $|V_i|^2= 3$
and  that 
the  complex valued vector 
$$U= \frac{z_1}{z_1-z_2}V_1+\frac{z_2}{z_2-z_3}V_2+\frac{z_3}{z_3-z_1}V_3,$$ has purely imaginary components: thanks to the fact that
$V_1+V_2+V_3=0$, and using $\bar  z_i= z_i^{-1}$, one sees that $\bar U= -U$. Therefore $-U\otimes U  \geq 0$. 

On the other hand, the tensor
$\frac{z_1z_2}{(z_1-z_2)^2}V_1^{\odot^{2}}+\frac{z_2z_3}{(z_2-z_3)^2}V_2^{\odot^{2}}+\frac{z_1z_3}{(z_3-z_1)^2}V_3^{\odot^{2}}$ is real.

\ermq
\bpf (Of Corollary~\ref{cor.eqwder}).
Equation~\eqref{eqwder1} is a direct consequence of Lemma~\ref{lem.derz}, while~\eqref{eqwder2} follows from
\beqnas\nabla\sigma\otimes \nabla \sigma=&& -[\frac{z_1E_1-z_2E_2}{z_1-z_2}+\frac{z_2E_2-z_3E_3}{z_2-z_3}+\frac {z_3E_3-z_1E_1}{z_3-z_1}]^{\otimes^{2}}\\
=&& -[E_1+\frac{z_2}{z_1-z_2}(E_1-E_2)+E_2+\frac{z_3}{z_2-z_3}(E_2-E_3)+E_3\\
&&+\frac{z_1}{z_3-z_1}(E_3-E_1)]^{\otimes^{2}}
\eeqnas 
and the fact that $E_1+E_2+E_3=0$.
\epf

First, observe that when $z_2= e^{i\phi} z_1$ with $\phi\to 0$,  with $z_1\neq 1,j, \bar j$, the tensor $\nabla\nabla \sigma$ is  equivalent to  $-\frac{4}{\sin^2\phi/2} V_1\otimes V_1$, and therefore is not bounded below. This shows that there cannot exist any $CD(\rho, \infty)$ inequality for $\lambda<1$.

From now on, the parameter  $a\in \bbR$ being fixed, let us call $b(a)$ the best constant $b$ in inequality~\eqref{ineg.gamma2.local} at a given point in the interior of the triangle. It  is obtained the lowest eigenvalue of the symmetric tensor
$\nabla\nabla\sigma- a\nabla\sigma\otimes\nabla \sigma$.  We want to understand when this function is bounded below.

It is then better to change coordinates and consider a basis $W_1= V_1/\sqrt{3}$ and $W_2$ which is orthogonal to $W_1$ and norm $1$, such that
$$V_2=\sqrt{3}(-\frac{1}{2} W_1+\frac{\sqrt{3}}{2}W_2), V_3=\sqrt{3}(-\frac{1}{2} W_1-\frac{\sqrt{3}}{2}W_2).$$

Moreover, we shall set $z_2= z_1u$, $|u|=1$, so that $z_3= 1/(z_1^2u)$, and $z_1^3=z$.  Observe that  the image of  $(x,y)\mapsto (z,u)$,  is $S_1^2$,  where $S_1$ is the set complex numbers with modulus 1.

With those notations, inequality~\eqref{ineg.gamma2.local} becomes 
$$\frac{AW_1^{\odot^2}+B W_2^{\odot^2}+C W_1\odot W_2}{4(1-u)^2(1-zu)^2(1-zu^2)^2} \geq b(W_1^{\odot^2}+W_2^{\odot^2}), $$

where 
 \beqnas
 A&=& 3\Big( \big(a(u^2+1)+(2a-4)u\big)(1+u^6z^4)  -uz(u+1)\big(4a(u^2+1)+u^2-10u+1 \big)(1+z^2u^3) \\&&
 +2u^2z^2\big( 2a(u^4+1)+a(u^3+u) +(6a-12)u^2  \big)\Big)
 \eeqnas
\beqnas
B&=&9(u-1)^2\Big(  a (u^6 z^4+1)-(u^5z^3+uz)-(u^4z^3+u^2z)+(-2a+4)u^3z^2 \Big)
 \eeqnas

\beqnas
C&=&6 \sqrt{3}(u-1)(1-z^2u^3)\big( a(u^4z^2+1)+a(u^3z^2+u)+(1-2a)(u^3z+uz)-2u^2z \big)
 \eeqnas
Setting $z= e^{i\theta}$ and $u= e^{i\phi}$ we have 
\beqnas
 A&=&12u^4z^2\Big(2\cos(2\theta+3\phi)(a\cos^2 \phi/2 -1)  -2\cos(\phi/2)\cos(\theta+3/2\phi)((4a+1)\cos \phi -5) \\&&
 +2a\cos(2\phi)+a\cos(\phi)+3(a-2)\Big)
\eeqnas

\beqnas
B&=&-72u^4z^2(\sin^2(\phi/2))\Big( a\cos(2\theta+3\phi)-\cos(\theta+2\phi) -\cos(\theta+\phi)+(2-a)  \Big)
\eeqnas

\beqnas
C&=&48\sqrt{3}z^2u^4 \sin(\phi/2)\sin(\theta+3/2\phi)\Big(a\cos(\theta+2\phi)+a\cos(\theta+\phi)+(1-2a)\cos(\phi)-1\Big)
\eeqnas

while the denominator may be written as
$$-4^4z^2u^4\sin^2(\phi/2)\sin^2(\theta/2+\phi/2)\sin^2(\theta/2+\phi).$$
Simplfying everything by $u^4z^2$, and letting $A_1= -Az^{-2}u^{-4}$, $B_1= -Bz^{-2}u^{-4}$, $C_1=Cz^{-2}u^{-4}$, and $N= 4^4\sin^2(\phi/2)\sin^2(\theta/2+\phi/2)\sin^2(\theta/2+\phi)$, we see that the best
constant $b(a)$, at some point $(z,u)$  is then 
$$b(a)= \frac{A_1+B_1-\sqrt{(A_1-B_1)^2+ C_1^2}}{2N}.$$

We now may prove the following

\blem\label{bnd.on.b}

The function $b(a)$ is unbounded below on the set $|z|=|u|=1$ if  $a>1/3$.

\elem

\bpf (Of Lemma~\ref{bnd.on.b})
We shall see  in Lemma~\ref{bnd.on.b.1/3} that the function $b(a)$ is bounded below for $a=1/3$. This of course shows that it is also bounded below for any $a<1/3$. 
What we have to prove then is that 
 the function $b(a)$ is unbounded below when $a>1/3$. 
 
 Let us concentrate on the case $a<1$. It is enough to observe the asymptotics of $b(a)$ around $\theta = \phi= 0$.  The result is obtained   when choosing  $\phi= \lambda \theta^2$.  Then, one has  $$\bcas
 A_1\simeq 12(1-a)\theta^4= \alpha\theta^4\\
 B_1\simeq 18\lambda^2(1-2a)\theta^6= \beta\theta^6\\
 C_1\simeq -24\sqrt{3} \lambda a\theta^5= \gamma\theta^5\ecas
 $$
 
Then,  $b(a)\simeq c/\theta^2 $, where the constant $c$ has the sign of $4\beta\alpha-\gamma^2$, that is of $(1-3a)$.  When $a>1/3$, this converges to $-\infty$ when $\theta\to 0$.

 \epf 
 
 \brmq When choosing in the previous argument $\phi= \lambda\theta$, one sees that $b(a)$ is unbounded below as soon as $a>1/2$.  This asymptotics is not enough to capture the optimal bound.

 \ermq

We now concentrate on the case $a= 1/3$. We are able to compute explicitly the lower bound for $a=1/3$, which corresponds and fits with the $SU(3)$ computation, although the explicit computation of the lower bound is not explicit (and not really of interest) for the other values of $a<1/3$.

We will study the function  $b(a)$ in case $a=\frac{1}{3}$.

\blem\label{bnd.on.b.1/3}
 The function $b(1/3)$ is bounded below and it's lower bound is $9/8$
 \elem

\bpf (Of Lemma~\ref{bnd.on.b.1/3})

In the case $a=1/3$ the function $b(1/3)$ have the following form in (z,u): 

\beqnas
b(1/3)&=&\frac{1}{2} \frac{1}{(u-1)^2(zu^2-1)^2(zu-1)^2}\Big(  P(z,u)-\sqrt{Q(z,u)}  \Big)
\eeqnas 
 where 
 \beqnas
 P(z,u)&=&(u^2-4u+1)(1+u^6z^4)-4zu(u+1)(u^2-3u+1)(1+z^2u^3)\\&&
 +u^2z^2(u^4+8u^3-30u^2+8u+1).\\
 Q(z,u)&=&[(z^2u^4-zu^3-zu^2+u^2-u+1)(z^2u^2-z^2u^3+z^2u^4-zu-zu^2+1)\\&&
(z^2u^4+z^2u^3+zu^3-6zu^2+zu+u+1)^2].
 \eeqnas

Using the same notations than in Lemma~\ref{bnd.on.b},  $b(1/3)$ can be written  also in $(\theta,\phi)$ in  this form 
\beq \label{ega} 
b(1/3)=\frac{1}{2}\frac {N(\theta, \phi)}{D(\theta, \phi)}
\eeq

where 
\beqnas
N(\theta, \phi)&=& 2\Big( 2\cos(2\theta +3\phi)(\cos(\phi)-2)-8\cos(\phi/2)\cos(\theta +3/2\phi)(2\cos(\phi)-3)\\&&
+\cos(2\phi)  +8\cos(\phi) -15 -\mid T \mid  \mid 2\cos(\phi/2)  \cos(\theta +3/2\phi)  +\cos(\phi)  -3 \mid \Big)\\
D(\theta, \phi)&=&-2^6 \sin^2(\phi/2)\sin^2(\theta/2+\phi)\sin^2(\theta/2+\phi/2)
\eeqnas 
where    $$T=z^2u^3-2\cos(\frac{\phi}{2})zu^{\frac{3}{2}}+2\cos(\phi)-1$$
and $\mid T \mid =\sqrt{ T \bar{ T}} $\\

  Setting $x=\cos(\phi/2), ~ y=\cos(\theta+3/2\phi )$, and $y=x+w$ we may rewrite this as 
 
 \beq\label{form}
 b(1/3)=\frac{1}{4}\frac{ (2(1-x^2)-xw)^2+3w^2(1-x^2)-(2(1-x^2)-xw)\sqrt{(2(1-x^2)-xw)^2+3w^2(x^2-1)}}{(1-x^2)w^2}
 \eeq
  
 To see this, we just replacing in \eqref{ega} 
 \beqnas
 &&\cos(\phi)=2x^2-1, ~\cos(2\phi)=  8x^4-8x^2+1, \cos(2\theta+3\phi)=2y^2-1 \\
&&\sin^2(\frac{\phi}{2})=1-x^2, \sin(\phi+\frac{\theta}{2})\sin(\frac{\theta+\phi}{2})=\frac{1}{4}(x-y)
\eeqnas
we obtain 
\beqnas
N(x,y)&=&8\Big(2x^4-8x^3y+2x^2y^2+x^2+10xy-3y^2-4 -\\
&&\mid x^2+xy-2 \mid  
\sqrt{4x^4+4x^2y^2-4x^3y-3y^2-7x^2+2xy+4} \Big)
\eeqnas
\beqnas
D(x,y)&=& -2^4(1-x^2)(x-y)^2
\eeqnas 
then if we set $y=x+w$  we have the result

 Now, in formula \eqref{form}, we set $t=\mid\frac{2(1-x^2)-xw)}{w\sqrt{1-x^2}}\mid$, and then $b(1/3)$ becomes  
 $$  b(1/3)=\frac{1}{4}(t^2+3-t\sqrt{t^2-3}), ~~t  \ge \sqrt{3}.$$ 
It is an easy exercise to check that the lower bound of this last function of $t$ is  $ 9/8$. 

\epf

We now collect the results of Lemmas~\ref{bnd.on.b} and~\ref{bnd.on.b.1/3} to get Proposition~ \ref{CD.r.n.brute.force}.
 
\epf
\subsection{ A simpler proof of the curvature-dimension inequality\label{subsec.simpler.proof}}

As mentioned in the introduction, we shall show that the use of the complex coordinates $(Z, \bar Z)$ provide a much simpler proof of Proposition~\ref{CD.r.n.brute.force}. Everything relies on the boundary equation~\eqref{eq.boundary}, which is nothing else than a particular case of a general equation which is valid as soon as orthogonal polynomials come into play (see~\cite{BOZ}). In the  coordinates $(Z, \bar Z)\in \cD$, the operator $\LL ^{(\lambda)}$ takes a simpler form, even if the metric looks  more complicated. This illustrates the use of the appropriate coordinates whenever one has a polynomial structure such as this deltoid model. 

Once again, our aim is to compute the Hessian of the function  $\log P$, where $P$ is defined in equation~\eqref{defP}.

Following~\cite{bglbook}, the Hessian of $f$, applied to $dh,dk$, that is in a local system of coordinates 
$\nabla\nabla^{ij}(f)\partial_ih\partial_j k$, may be defined as 
\beq\label{def.hessian}H[f](h,k)= \frac{1}{2}\big(\Gamma(h, \Gamma(f,k))+ \Gamma(k, \Gamma(f,h))- \Gamma(f,\Gamma(k,h)\big).\eeq
We want to apply this with $f = \log P$ and $h,k= Z,\bar Z$. For this, one may  use the boundary equation which takes in this context the particular form
\beq\label{eq.boundary}\Gamma(Z,\log (P))= -3Z, ~\Gamma(\bZ, \log P)= -3\bZ,\eeq and is easily checked from formulae~\eqref{eq.deltoid}

From this, we deduce that, for any function $G(Z,\bZ)$, $\Gamma(\log P, G)= -3 D(G)$, where 
$D$ is the Euler operator $Z\partial_Z+ \bZ\partial_{\bZ}$

Let us write 
$$H^{11}=H[\log P](Z,Z),~ H^{12}= H[\log P](Z, \bZ), ~H^{22}= H[\log P](\bZ,\bZ).$$

From the previous remarks, we get
$$H^{11}= -3\Gamma(Z,Z) + \frac{3}{2}D(\Gamma(Z,Z)).$$
$$H^{12}= -3\Gamma(Z,\bZ) + \frac{3}{2} D(\Gamma(Z,\bZ)),$$ and 
$$H^{22} = -3\Gamma(\bZ,\bZ)+ \frac{3}{2}D(\Gamma(Z,Z)).$$
In other words, with the obvious notations,
$H= -3 \Gamma+ \frac{3}{2} D\Gamma$

In the same way, the tensor $\nabla\log P\otimes \nabla\log P$ may be computed in this system of coordinates as 
$$M= 9 \bpm Z^2& Z\bZ\\ Z\bZ& \bZ^2\epm$$.

In the end, the inequality
\beq\label{ineq.a1b1}-\nabla\nabla \log P \geq b_1 \Gamma + a_1 \nabla\log P\otimes \nabla \log P\eeq
amounts to 
$$(3-b_1)\Gamma - \frac{3}{2} D\Gamma - 9 a_1M\geq 0,$$
where $\Gamma$ denotes the matrix 
$$\bpm \Gamma(Z,Z)&\Gamma(Z,\bZ)\\ \Gamma(Z,\bZ)& \Gamma(\bZ,\bZ)\epm.$$
For such a tensor $(R^{ij})$ in complex coordinates, to represent a non negative real tensor amounts to ask that 
$$R^{12}\geq 0 \hbox{ and }
(R^{12})^2\geq R^{11}R^{22}.$$

$R^{12}\geq 0$ reads $3-b_1+ (b_1-18a_1)Z\bZ\geq 0$, and for this to be true on $\Omega$ amounts to ask 
$$b_1\leq 3, \quad a_1\leq 1/6,$$ since $Z\bZ$ varies from $0$ to $1$ on $\Omega$.

The second one writes
\beq\label{cond.cdZbZ}[(3-b_1)/2+ b_1/2-9a_1 Z\bZ]^2\geq (3/2-b_1)^2Z\bZ+ (b_1-9a_1)^2Z^2\bZ^2+ (Z^3+\bZ^3)(b_1-9a_1)(3/2-b_1).\eeq
Writing everything in polar coordinates $Z= \rho e^{i\theta}$, this writes as 
$$P_2(\rho ^2)\geq 2 \rho^3\cos(3\theta)(b_1-9a_1)(3/2-b_1),$$
where $P_2$ is a degree 2 polynomial.

Observe that this requires to be true for any $\rho \in [0,1]$ when $\cos(3\theta)=0$ (which corresponds to the cusps of the deltoid curve).

But, with the explicit computation of $P_2$, one gets 
$$P_2(\rho^2)-2 \rho^3(b-9a)(3/2-b)=\frac{1}{4}(1-\rho)(3-b_1+ b_1\rho)\big(3-b_1+ \rho(3-2b_1)+ 3\rho^2(b_1-12a_1)\big).$$

For the maximal value $a_1= 1/6$, 
$$P_2(\rho^2)-2 \rho^3(b-9a)(3/2-b)= (1-\rho)^2(3-b_1+ 3\rho(2-b_1)),$$ and we get a bound $b_1\leq 9/4$.

For these values $a_1=1/6$ and $b_1= 9/4$, equation~\eqref{cond.cdZbZ}  writes
$$4(1-\rho^2)^2 \geq \rho^2 +\rho^4-2\rho^3\cos(3\theta),$$
while the condition $\Gamma(Z,\bZ)^2\geq \Gamma(Z,Z)\Gamma(\bZ,\bZ)$, which characterizes the points in $\bar \cD$, writes
$$\frac{1}{4}(1-\rho^2)^2 \geq \rho^2 +\rho^4-2\rho^3\cos(3\theta),$$ so the the inequality is satisfied everywhere in $\cD$. Observe that the critical points for the curvature-dimension inequality for the critical values are attained at the cusps.

\brmq
 Observe that the values $b_1= 9/4$ and $a_1= 1/6$ in inequality~\eqref{ineq.a1b1} are once again exactly the bounds obtained in equation~\eqref{ineq.gamma2.w}. Moreover, we know that even with $a_1=0$ (corresponding to a $CD(\rho, \infty)$ inequality), if we look for the optimal value for $b_1$, it is clear from this method that the best constant $b_1$ is $b_{max}< 3$, so that whatever the constant $a_1\in [0, 1/6]$, the optimal value for $b_1$ lies in the interval $[9/4, 3)$. (The optimal constant $b_1(a_1)$ may be explicitly computed but has no real interest.)
\ermq

\section{Sobolev inequalities and bounds on the eigenvectors\label{sec.bound.on.polyn}}

As described in Section~\ref{sec.cdrn}, from the curvature dimension inequality, we may obtain bounds on the supremum of the associated eigenvectors. 
More precisely,  whenever a $CD(\rho,n)$ inequality holds with $\rho>0$ and $n<\infty$, there exists a constant $C$ such that  for any eigenvector $P$ satisfying $\LL (P)= -\mu P$, then  
$$\|P\|_\infty \leq  C \mu^{n/4}\|P\|_2,$$where the $\cL^ 2$ norm is computed with respect to the invariant measure of the operator $\LL $.

Turning to the case of the operator $\LL ^{(\lambda)}$ on the deltoid, we  recall from \cite{Zribi2013} that the associated eigenvectors  which are polynomials with total degree $n$ have eigenvalues $\mu_{p,q}= (\lambda-1)(p+q)+p^2+q^2+pq$, with $p+q=n$.  More precisely, of any $n\geq 1$ such that $p+q=n$,  when $p\neq q$, there is a dimension $2$ associated eigenspace. In complex variables, for such value $\mu_{p,q}$, there is a unique degree $n$ polynomial $P_{p,q}(Z, \bar Z)$ with highest degree term $Z^p\bar Z^q$ and another one which is $\bar P_{p,q} (Z,\bar Z)= P_{p,q} (\bar Z, Z)= P_{q,p}(Z, \bar Z)$ eigenvector (the polynomial $P_{p,q}$ having real coefficients).  For $p=q$ however, the 	associated eigenspace is one dimensional.  The real forms are $S_{p,q}= \frac{1}{2}(P_{p,q}+ P_{q,p})$ and $A_{p,q}= \frac{-i}{2}(P_{p,q}-P_{q,p})$, which  form a real basis for this eigenspace.

When $\lambda>1$, for $\LL ^{(\lambda)}$, for any $\mu_{p,q}$ and of any polynomial $P$ in the associated eigenspace, one gets from the $CD(\frac{3(\lambda-1)}{4}, 2\lambda)$ inequality 
\beq\label{bnd.eigenv1}\| P\|_\infty \leq C(\lambda) \mu^{\lambda/2}.\eeq

Looking at the constants, this does not produce  any estimates for $\lambda=1$ or $0<\lambda <1$. However, for $\lambda=1$, one may consider the following.
The operator $\LL ^{(1)}$ is nothing else than the usual Laplace operator acting on functions $f(Z)$, where the function $Z$ is given in~\eqref{def.zk}. As functions of $(x,y)$ in the real plane, those functions are periodic in $x$ with period $4\pi$ and in $y$ with period $4\pi/\sqrt{3}$. As such, the associated semigroup $P_t^{(1)}$  is an image of the 
product semigroup of the associated $1$ dimensional torus (that is the semigroup on the real line acting on periodic functions).  More precisely, when considering a function on the deltoid as a function of $(x,y)\in \bbR^2$, 
$$P_t^{(1)}((x,y, dx', dy')= P_t^{S^1(4\pi)}(x,dx')P_t^{S^1(4\pi/\sqrt{3})}(y, dy'),$$ where 
$P_t^{S^1(\tau)}(x,dx')$ is the semigroup of the torus with radius $\tau$, that is the semigroup of the one-dimensional Brownian motion acing on $\tau$-periodic functions. Since both semigroups have a density which is bounded above by $C/\sqrt{t}$ for some constant $C$ and  for $0<t\leq 1$, it turns out that $P_t^{(1)}$ has a density which is bounded above by $C'/t$. This is enough to get the bound on the associated eigenvectors. In the end, we get
\bprop  For any $\lambda\geq 1$, there exists a constant $C(\lambda)$ depending on $\lambda$ only,  such that for  any polynomial $P$ eigenvector of $\LL ^{(\lambda)}$ with eigenvalue $\mu\neq 0$, one has 
\beq\label{bnd.eigenvec}\| P\|_\infty \leq C(\lambda) \mu^{\lambda/2} \|f\|_2.\eeq

\eprop

\brmq Looking at the constants, whenever $\lambda\to 1$, the constant $C(\lambda)$ in~\eqref{bnd.eigenvec} goes to $\infty$, and there is    an unexpected discontinuity in the constants. Indeed, our computations are not the best possible. One may sharpen them with the help of spectral gaps, that is the knowledge of the lowest non $0$ eigenvalue,  which here is  $\lambda$. More precisely, one may reinforce the constants in the ultracontractive bounds  under a $CD(\rho, n)$ inequality and the knowledge of this lowest eigenvalue (see remark page 313 in~\cite{BGL}). But the argument in this modified estimate produces a Sobolev inequality with any dimensional parameter $m>n$ , and a constant which is not improved when $m\to n$.  There is therefore a balance in the optimal bound on $P$ between the value of $\mu$ (for $\mu$ large one wants $m$ to be the lowest possible), and for $\lambda\to 1$ (when $\lambda\to 1$, one wants $C(\lambda)$ not too big).   We could  such have produced a better bound. But indeed, the remark  in~\cite{BGL}  as it stands is not really valid for $\rho=0$ which corresponds in our case to $\lambda=1$. One would have to sharpen this estimate, both for the case $\rho=0$ and for the value of $m$. It is indeed true that one may obtain a $n$-dimensional Sobolev inequality   (under it's entropic form) under a estimate on the lowest eigenvalue and some $CD(\rho,n)$ inequality, even for $\rho<0$, but the argument in~\cite{BGL} is clearly not sufficient for that and requires further analysis.
\ermq

From the point of view of the invariant measure $\mu^{(\lambda)}$ of $\LL ^{(\lambda)}$, what is relevant is the decomposition of $\cL^2(\mu^{(\lambda)})$ into spaces of orthogonal polynomials. More precisely,  when denoting $\cP_k$ the space of polynomials with total degree less than or equal to $k$, one considers the subspace $\cH_k$ of $\cP_k$ which is orthogonal to $\cP_k$, such that one has the orthogonal decomposition
$$\cL^2\big(\mu^{(\lambda)}\big)= \oplus_{k=1}^\infty \cH_k,$$ where $\cH_0$ is the space of constant functions. 

One has
\bprop There exists a constant  $C_1(\lambda)$ such that, for any $k\geq 1$ and any $P\in \cH_k$
\beq\label{bnd.Hk}\| P\|_\infty \leq C_1(\lambda)k^{\lambda +1/2}\|P\|_2.\eeq

\eprop
 
 \bpf One may decompose $\cH_k$ into the eigenspaces associated to $\LL ^{(\lambda)}$. There are  $r_k=[k/2]+1$ such eigenspaces, and all the eigenvalues belong to the interval $[k(\lambda+k-5/4), k(\lambda+k-1)]$, or, when $k\geq 1$, in the interval $[3/4k^2, \lambda k^2]$.

 Writing $P\in \cH_k$ as 
 $P= \sum_{i=1}^{r_k} a_iP_i$ where $P_i$ is an eigenvector with $\|P_i\|_2=1$ and $\|P\|^2= \sum_1^{r_k} a_i^2$, one has from~\eqref{bnd.eigenvec} and the bound on the eigenvalues in $\cH_k$
 \beqnas \|P\|_\infty &&\leq \sum_1^{r_k} |a_i| \|P_i\|_\infty \leq  C(\lambda) \lambda^{(1+\lambda/2)} \sum_1^{r_k} |a_i| k^{\lambda}\\&& \leq C(\lambda)\lambda^{(1+\lambda/2)}k^{\lambda}r_k^{1/2}\sqrt{\sum_1^{r_k} a_i^2},
 \eeqnas
 from which the bound follows immediately.

 \epf
 
 One may wonder how far inequality~\eqref{bnd.Hk} is from the Sobolev inequality we started from. Observe first that, for the heat kernel $\cP_t^{(\lambda)}$, for any function $P\in \cH_k$, one has 
 $$\|\cP_t^{(\lambda)} P\|_2 \leq \exp(-3/4t k^2) \|P\|_2,$$ since all eigenvalues of $\cP_t$ on $\cH_k$ are bounded below by $\exp(-3/4 tk^2)$.
 
 Therefore, we have, for any $P\in \cH_k$
 $$\|\cP_t^{(\lambda)} P\|_\infty \leq C_1(\lambda) \exp(-3/4t k^2) k^{\lambda+1/2}\|P\|_2.$$
 Observe that this relies only on the bound~\eqref{bnd.eigenvec} together with the knowledge of the eigenvalues.

\bthm\label{recque.sob} Let  $P_t$ be a symmetric Markov semigroup  with reversible probability measure $\mu$ and generator $\LL$. Assume that  $\LL$ satisfies a Poincaré inequality  and that  one has a decomposition  into orthogonal spaces 
$\cL^2(\mu) = \oplus_k \cH_k$, where $\cH_k$ is a linear space, with the property that,  for some real number $a>0$  and  for any 
$f\in \cH_k$,  
$$\| P_t f\|_\infty \leq Ck^pe^{-atk^2} \|f\|_2.$$
Then,
$\LL $ satisfies a tight Sobolev inequality with dimension $m=2p+1$.

\ethm 

\bpf Following~\cite{BGL}, and from the existence of a Poincaré inequality, it is enough to prove that, for $t\in (0,1]$ and  for some constant $C$,
$\| P_t f\|_\infty  \leq C t^{-m/4} \|f\|_2$. We may  restrict to the case where $\|f\|_2=1$. 
For $f\in \cL^2(\mu)$, let us write 
$f= \sum_k f_k$, where $f_k\in \cH_k$ and  $\sum_k \|f_k\|_2^2= 1$.  
$$\|P_t f\|_\infty \leq \sum_k \|P_t f_k\|_\infty \leq \sum_k k^p e^{-at k^2} \|f_k\|_2 \leq \Big( \sum_k k^{2p} e^{-2atk^2}\Big)^{1/2}.$$

One may compare  the sum  $\sum_k k^{2p} e^{-2atk^2}$ with $\int_0^\infty x^{2p}\exp(-2atx^2) \, dx$, where the function $x^{2p}\exp(-2atx^2)$ is increasing on $(0, \sqrt{p/(2at)}$ and decreasing on $(\sqrt{p/(2at)}, \infty)$, and we see that, for $0<t\leq 1$, 
$$\sum_k k^{2p} e^{-2atk^2}\leq C(a,p) t^{-(p+1)/2}, ~ 0<t \leq 1.$$
Therefore, following the results exposed in Section~\ref{sec.cdrn}, we get a Sobolev inequality with dimension $m= 2(p+1)$. The existence of  a Poincaré inequality (that is of a strictly positive first non zero eigenvalue for $-\LL $) insures that we may get a tight Sobolev inequality~\eqref{def.sob}.
 This gives the result.

\epf

Looking at the values for the deltoid model, we see that the estimate provides a Sobolev inequality with dimension $m= 2\lambda+3$, whereas we started from a Sobolev inequality with dimension $2\lambda$. One may wonder if this lost in dimension (from $n$ to $n+3$) is due to too crude estimates on both the eigenvalues and the summation in the series, or from the fact that the spaces $\cH_k$ are $k$ dimensional. Indeed, even in the case or one dimensional Jacobi operators, where the eigenspaces are one dimensional, where  the eigenvalues for the associated operator are $k(k+c)$ for polynomials with degree $k$, one would pass with the same method  from a Sobolev inequality with dimension  $n$ to a Sobolev  inequality with dimension $n+1$. 
  This is in big contrast with the case of logarithmic Sobolev inequalities, where  estimates on the $\cL^p$ bounds on the eigenvectors are indeed equivalent to  logarithmic Sobolev inequalities (see~\cite{bglbook}).

Finally, we directly get from this  a criterium for a symmetric operator constructed from orthogonal polynomial would have a bounded density.

\bprop
Let $K$ be a symmetric operator in $\cL^2(\mu^{(\lambda)})$ which maps $\cH_k$ into $\cH_k$ and is such that, for any $P\in \cH_k$, $
\| K(P)\|_2 \leq \nu_k \|P\|_2$. If 
$A= \sum_k \nu_k^2 k^{2\lambda+1} < \infty$, then 
$K^2$ may be represented by a bounded kernel, 
$$K^2(f)(x)= \int f(y) k(x,y) \, d\mu(y),$$
where $|k| \leq A$.

\eprop

\bpf
Arguing as in the proof of Theorem~\eqref{recque.sob}, we may write $f= \sum_k f_k$ with $f_k\in \cH_k$. Then,
\beqnas \|K(f)\|_\infty &&\leq \sum_k \| K(f_k)\|_\infty\leq \sum_k k^{\lambda/1/2}\| K(f)\|_2\\&&\leq \sum_k \nu_kk^{\lambda/1/2}\|f_k\|_2\leq (\sum_k \nu_kk^{\lambda/1/2})^{1/2} \|f\|_2.
\eeqnas
 Therefore, the operator $K$ is bounded from  $\cL^2$ into $\cL^\infty$ with norm $A^{1/2}$. By symmetry and duality, the same is true from $\cL^1$ into $\cL^2$, and by composition, $K^2$ is bounded from $\cL^1$ into $\cL^\infty$ with norm $A$. It therefore may be represented by a kernel $k$ bounded by $A$. 

\epf

\brmq\label{rmq.nash} The method presented here says nothing about the case where $0<\lambda< 1$. Indeed, in this case, one may expect to have a two-dimensional behavior for the heat kernel, that is $\|P_1\|_{2, \infty} \leq Ct^{-1/2}$, $0<t\leq 1$. In this context, it is better to replace Sobolev inequalities by Nash inequalities, that is inequalities of the form
$$\|f\|_2^2 \leq \|f\|_1^{2\theta}\Big(\|f\|_2^2+ C \int \Gamma(f,f)d\mu\Big)^{1-\theta},$$ where $\theta= \frac{2}{n+2}$ is a dimensional parameter. When $n>2$, this is equivalent to a Sobolev inequality with dimension $n$, but for $n\in (1,2)$, this is still equivalent to a bound $\|P_t\|_{2,\infty}\leq C't^{-n/4}$ (see~\cite{bglbook}). As mentioned, we may expect when $0<\lambda< 1$ some Nash inequality with dimensional parameter $n=2$. We cannot expect any smaller value for $n$ since, applied to any function compactly supported in the interior of $\cD$, this would contradict the classical two dimensional Nash inequality in an open domain of $\bbR^2$.  However, the singularity of the measure at the cusps of the deltoid make things a bit hard to analyze. Following the method developed in~\cite{bglbook},  pages 370-371, we are able to prove Nash inequalities with dimension $n= 5/2$, with however  a constant $C(\lambda)$  which goes to infinity when $\lambda\to 0$. This is not satisfactory for many points of view. First, because of the bad dimension, and for the lack of continuity in these inequalities when $\lambda\to1$.  Secondly, when $\lambda\to 0$, the measure converges to  the uniform measure on the three cusps $\frac{1}{3}(\delta_{1/3}+\delta_{j/3}+ \delta_{\bar j/3})$,  and the associated polynomials converge to the corresponding polynomials on three point (when one would have to replace $Z\bar Z$ by $1/9$ and $Z^3$ by $Z/3^3$). It is therefore a challenging question to produce for small $\lambda$ some functional inequality which recaptures this particular structure in the limit. 

\ermq

 \bibliographystyle{amsplain}   
\bibliography{Bib.these.OZ.cd}
\end{document}